\pdfoutput=1

\documentclass[proceedings]{aofa}

\usepackage{tikz}
\usetikzlibrary{arrows,backgrounds,calc}
\usepackage{calrsfs}
\usepackage{hyperref}
\usepackage{lineno}
\usepackage{enumerate}
\usepackage{mathtools}
\usepackage{pgfplots}
\pgfplotsset{compat=1.5}

\makeatletter
\let\DMTCSproof\proof
\let\DMTCSendproof\endproof
\let\proof\@undefined
\let\endproof\@undefined
\makeatother

\usepackage{amsmath, amsfonts, amssymb, amsthm, amstext, bbm}

\let\proof\DMTCSproof
\let\endproof\DMTCSendproof

\allowdisplaybreaks

\newtheorem{theorem}{Theorem}

\newtheorem{proposition}{Proposition}[section]
\newtheorem{corollary}[proposition]{Corollary}

\theoremstyle{remark}
\newtheorem*{remark}{Remark}

\theoremstyle{definition}
\newtheorem*{definition}{Definition}

\DeclareMathOperator{\Reg}{Reg}

\newcommand{\FGF}{H}
\newcommand{\TODO}[1]%
{\par\fbox{\begin{minipage}{0.9\linewidth}\textbf{TODO:} #1\end{minipage}}\par}

\renewcommand{\P}[0]{\mathbb{P}}
\newcommand{\E}[0]{\mathbb{E}}
\newcommand{\V}[0]{\mathbb{V}}
\newcommand{\N}[0]{\mathbb{N}}

\DeclareMathOperator{\cdeg}{cdeg}
\DeclarePairedDelimiter{\abs}{\lvert}{\rvert}

\newcommand{\upa}[0]{\mathnormal\uparrow}
\newcommand{\righta}[0]{\mathnormal\rightarrow}
\newcommand{\downa}[0]{\mathnormal\downarrow}
\newcommand{\lefta}[0]{\mathnormal\leftarrow}

\newcommand{\convexpath}[2]{
[
    create hullnodes/.code={
        \global\edef\namelist{#1}
        \foreach [count=\counter] \nodename in \namelist {
            \global\edef\numberofnodes{\counter}
            \node at (\nodename) [draw=none,name=hullnode\counter] {};
        }
        \node at (hullnode\numberofnodes) [name=hullnode0,draw=none] {};
        \pgfmathtruncatemacro\lastnumber{\numberofnodes+1}
        \node at (hullnode1) [name=hullnode\lastnumber,draw=none] {};
    },
    create hullnodes
]
($(hullnode1)!#2!-90:(hullnode0)$)
\foreach [
    evaluate=\currentnode as \previousnode using \currentnode-1,
    evaluate=\currentnode as \nextnode using \currentnode+1
    ] \currentnode in {1,...,\numberofnodes} {
-- ($(hullnode\currentnode)!#2!-90:(hullnode\previousnode)$)
  let \p1 = ($(hullnode\currentnode)!#2!-90:(hullnode\previousnode) - (hullnode\currentnode)$),
    \n1 = {atan2(\y1, \x1)},
    \p2 = ($(hullnode\currentnode)!#2!90:(hullnode\nextnode) - (hullnode\currentnode)$),
    \n2 = {atan2(\y2, \x2)},
    \n{delta} = {-Mod(\n1-\n2,360)}
  in
    {arc [start angle=\n1, delta angle=\n{delta}, radius=#2]}
}
-- cycle
}

\title[Reductions of Binary Trees and Lattice Paths]
{The Register Function and Reductions of\\ Binary Trees and Lattice Paths}

\author[B.~Hackl, C.~Heuberger, H.~Prodinger]{Benjamin Hackl\addressmark{1}\addressmark{\dag}
  \and Clemens Heuberger\addressmark{1}\thanks{B.~Hackl and C.~Heuberger are supported by
    the Austrian Science Fund (FWF): P~24644-N26 and by the Karl Popper Kolleg
    ``Modeling-Simulation-Optimization'' funded by the Alpen-Adria-Universit\"at Klagenfurt
    and by the Carinthian Economic Promotion Fund (KWF).}
  \and Helmut Prodinger\addressmark{2}\thanks{H.~Prodinger is supported by an incentive
    grant of the National Research Foundation of South Africa.}}

\address{\addressmark{1}Institut f\"ur Mathematik,
  Alpen-Adria-Uni\-ver\-si\-t\"at Klagenfurt, Austria, \{benjamin.hackl, clemens.heuberger\}@aau.at\\
  \addressmark{2}Department of Mathematical
  Sciences, Stellenbosch University, South Africa, hproding@sun.ac.za
}


\begin{document}

\maketitle
\begin{abstract}
  \paragraph{Abstract.}
  The register function (or Horton-Strahler number) of a binary tree is a well-known combinatorial parameter. We
  study a reduction procedure for binary trees which offers a new interpretation for the
  register function as the maximal number of reductions that can be applied to a given
  tree. In particular, the precise asymptotic behavior of the number of certain
  substructures (``branches'') that occur when reducing a tree repeatedly is determined.

  In the same manner we introduce a reduction for simple two-dimensional lattice
  paths from which a complexity measure similar to the register function can be
  derived. We analyze this quantity, as well as the (cumulative) size of an (iteratively)
  reduced lattice path asymptotically.
\end{abstract}

\section{Introduction}
\label{sec:introduction}

Binary trees are either a leaf or a root together with a left and a right subtree which
are binary trees. It is well-known that the generating function counting these objects
with respect to the number of inner nodes is given by
\[ B(z) = \frac{1 - \sqrt{1 - 4z}}{2z} = \sum_{n\geq 0} \frac{1}{n+1} \binom{2n}{n}
  z^{n}. \]
Thus, the $n$th Catalan number $C_{n} = \frac{1}{n+1} \binom{2n}{n}$ counts the number of
binary trees with $n$ inner nodes.

By simple algebraic manipulations, it is easy to verify that
$B(z)$ fulfills the identity
\[ B(z) = 1 + \frac{z}{1 - 2z} B\Big(\frac{z^{2}}{(1 - 2z)^{2}}\Big).  \]
However, as we will see in Section~\ref{sec:reduction-register}, we
can justify this identity from a combinatorial point of view as well, and the most
important part of this combinatorial interpretation is a reduction procedure for 
binary trees.

The aim of this paper is to analyze the binary tree reduction with a focus on the
structures that emerge when repeatedly reducing a given tree. After the aforementioned
introduction of the reduction in Section~\ref{sec:reduction-register}, we discover an
inherent connection to a very well-known branching complexity measure of binary trees: the
register function. 

Sections~\ref{sec:r-branches} and \ref{sec:all-branches} deal with the analysis of the
number of $r$-branches and the number of all branches within trees of given size, where an
$r$-branch can be thought of a local structure in a binary tree that survives exactly $r$
reductions.

In Section~\ref{sec:paths}, we switch our attention from binary trees to two-dimensional
lattice paths. As we will see, the generating function of these objects fulfills a similar
functional equation as the generating function for binary trees---and its combinatorial
interpretation strongly depends on a reduction process as well. The remainder of
Section~\ref{sec:paths} is devoted to analyzing the lattice path reduction. In particular,
Section~\ref{sec:fringes} investigates \emph{fringes} of lattice paths, which play a similar role
as branches with respect to binary trees.

On a general note, we used the open-source mathematics software system
SageMath~\cite{sage-7.0} in order
to perform the computationally intensive parts of the asymptotic analysis for each of the
quantities investigated in this paper. Furthermore, the proofs and many details are
omitted in this extended abstract; they can be found in the full version.

\section{Tree Reductions and the Register Function}
\label{sec:reduction-register}

As mentioned in the introduction, we want to find a combinatorial proof for the
following proposition.

\begin{proposition}\label{prop:tree-identity}
  The generating function counting binary trees by the number of inner nodes, $B(z) =
  \frac{1 - \sqrt{1 - 4z}}{2z}$, fulfills the identity
  \begin{equation}\label{eq:u1}
    B(z)=1+\frac{z}{1-2z}B\Big(\frac{z^2}{(1-2z)^2}\Big).
  \end{equation}
\end{proposition}

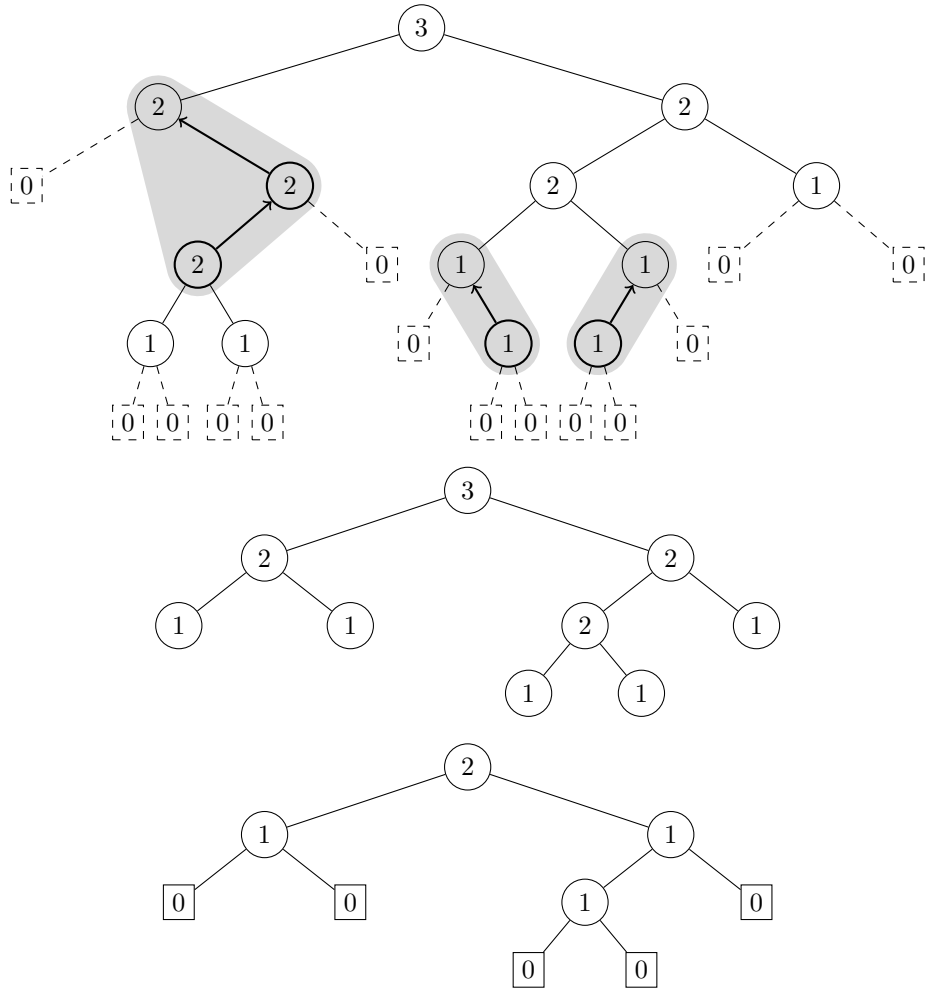
\begin{figure}[htbp]
  \centering
  \begin{tikzpicture}[scale=0.7,level distance=15mm,
    level 1/.style={sibling distance=100mm},
    level 2/.style={sibling distance=50mm},
    level 3/.style={sibling distance=35mm},
    level 4/.style={sibling distance=18mm},
    level 5/.style={sibling distance=8.5mm},
    ]
    \node[circle,draw] {${3}$}
    child {node[circle,draw] (n1) {${2}$} 
      child[dashed] { node [rectangle,draw]{${0}$}}
      child[<-, thick] {node [circle,draw] (n2) {${2}$}
        child {node[circle,draw] (n3) {${2}$}
          child[-,thin] {node[circle,draw] {${1}$}
            child[dashed] {node[rectangle, draw] {${0}$}}
            child[dashed] {node[rectangle, draw] {${0}$}}
          }
          child[-,thin] {node[circle,draw] {${1}$}
            child[dashed] {node[rectangle, draw] {${0}$}}
            child[dashed] {node[rectangle, draw] {${0}$}}
          }
        }
        child[-,dashed,thin] {node [rectangle,draw]{${0}$}
        }
      }
    }
    child {node[circle,draw] {${2}$}
      child {node [circle,draw]{${2}$}
        child {node[circle,draw] (e11) {${1}$}
          child[dashed] {node[rectangle,draw] {${0}$}}
          child[<-, thick] {node [circle,draw] (e12) {${1}$}
            child[dashed, thin, -] {node[rectangle,draw] {${0}$}}
            child[dashed, thin, -] {node[rectangle,draw] {${0}$}}
          }
        }
        child {node [circle,draw] (e21) {${1}$}
          child[<-, thick] {node [circle,draw] (e22) {${1}$}
            child[dashed, -, thin] {node [rectangle,draw]{${0}$}}
            child[dashed, -, thin] {node [rectangle,draw]{${0}$}}
          }
          child[dashed] {node[rectangle,draw] {${0}$}}
        }
      }
      child {node [circle,draw]{${1}$}
        child[dashed] {node[rectangle,draw] {${0}$}}
        child[dashed] {node[rectangle,draw] {${0}$}}
      }
    };
    \begin{pgfonlayer}{background}
    \fill[gray, opacity=0.3] \convexpath{n1, n2, n3}{17pt};
    \fill[gray, opacity=0.3] \convexpath{e11, e12}{17pt};
    \fill[gray, opacity=0.3] \convexpath{e21, e22}{17pt};
    \end{pgfonlayer}
  \end{tikzpicture}

  \vspace{2ex}

  \begin{tikzpicture}[scale=0.6,level distance=15mm,
    level 1/.style={sibling distance=90mm},
    level 2/.style={sibling distance=38mm},
    level 3/.style={sibling distance=25mm},
    level 4/.style={dashed,sibling distance=11mm},
    ]
    \node[draw,circle] {${3}$}
    child {node[draw,circle] {${2}$} 
      child{ node[draw,circle] {${1}$}
      }
      child {node[draw,circle] {${1}$}}
    }
    child {node[draw,circle] {${2}$}
      child {node[draw,circle] {${2}$}
        child {node[draw,circle] {${1}$}
        }
        child {node[draw,circle] {${1}$}
        }
      }
      child {node[draw,circle] {${1}$}
      }
    };
  \end{tikzpicture}

  \vspace{2ex}

  \begin{tikzpicture}[scale=0.6,level distance=15mm,
    level 1/.style={sibling distance=90mm},
    level 2/.style={sibling distance=38mm},
    level 3/.style={sibling distance=25mm},
    level 4/.style={dashed,sibling distance=11mm},
    ]
    \node[draw,circle] {${2}$}
    child {node[draw,circle] {${1}$} 
      child{ node[draw,rectangle] {${0}$}
      }
      child {node[draw,rectangle] {${0}$}}
    }
    child {node[draw,circle] {${1}$}
      child {node[draw,circle] {${1}$}
        child {node[draw,rectangle] {${0}$}
        }
        child {node[draw,rectangle] {${0}$}
        }
      }
      child {node[draw,rectangle] {${0}$}
      }
    };  
  \end{tikzpicture}
  \caption{Illustration of the compactification $\Phi$: in the first tree, the leaves are
    deleted (dashed nodes) and nodes with exactly one child are merged (gray overlay). The
  second tree shows the result of these operations. Finally, in the last tree all nodes
  without children are marked as leaves.}
  \label{fig:trimming}
\end{figure}

\begin{proof}[(Sketch)]
  The central idea of this proof is to consider a reduction of a binary tree $t$, which we
  write as $\Phi(t)$:

  First, all leaves of $t$ are erased. Then, if a node has only one child, these two nodes
  are merged; this operation will be repeated as long as there are such nodes. Finally,
  the nodes without children are declared to be leaves.

  Observe that this reduction is only defined for trees $t$ that have at least one inner
  node. The various steps of this operation (which was introduced in \cite{YaYa10}) are
  depicted in Figure~\ref{fig:trimming}. The number attached to the nodes will be
  explained later. 

  It can be shown that the generating function
  \[ \frac{z}{1 - 2z} B\big(\frac{z^{2}}{(1 - 2z)^{2}}\big)  \]
  counts all binary trees that can be reduced at least once. Thus, the functional
  equation~\eqref{eq:u1} can be interpreted combinatorially as follows: a binary tree is
  either just $\square$, or it can be reduced at least once.
\end{proof}
\begin{remark}
  Note that \eqref{eq:u1} can be used to find a very simple proof for a well-known
  identity for Catalan numbers: 

  Comparing the coefficients of $z^{n+1}$, \eqref{eq:u1} leads to
  \begin{align*}
    C_{n+1}&=[z^{n+1}]\sum_{k\ge0}C_k\frac{z^{2k+1}}{(1-2z)^{2k+1}}=
             \sum_{k\ge0}C_k[z^{n-2k}]\sum_{j\ge0}2^{j}\binom{2k+j}{j}z^j\\
           &=\sum_{0\le k\le n/2}C_k2^{n-2k}\binom{n}{2k},
  \end{align*}
which is known as Touchard's identity~\cite{shapiro-1976, touchard-1924}.
\end{remark}

With this interpretation in mind, \eqref{eq:u1} can also be seen as a recursive process to
generate binary trees by repeated substitution of chains. This process can be modeled by
the generating functions
\begin{equation}\label{eq:register-leq-recursion}
B_0(z)=1, \quad B_r(z)=1+\frac{z}{1-2z}B_{r-1}\Big(\frac{z^2}{(1-2z)^2}\Big), \quad r\ge1.
\end{equation}
By construction, $B_{r}(z)$ is the generating function of all binary trees that can be
constructed from $\square$ with up to $r$ expansions---or, equivalently---all binary trees
that can be reduced to $\square$ by applying $\Phi$ up to $r$ times.

As it turns out, these generating functions are inherently linked with the \emph{register
  function} (also known as the Horton-Strahler number) of binary trees. In order to
understand this connection, we introduce the register function and prove a simple property
regarding the compactification $\Phi$.

\begin{figure}[htbp]
  \centering
  \begin{tikzpicture}[scale=0.72,level distance=15mm,
    level 1/.style={sibling distance=90mm},
    level 2/.style={sibling distance=40mm},
    level 3/.style={sibling distance=30mm},
    level 4/.style={sibling distance=15mm},
    level 5/.style={sibling distance=8mm},
    ]
    \node[circle,draw,color=magenta] {${3}$}
    child {node[circle,draw,color=green] {${2}$} 
      child[color=black]{ node [rectangle,draw,color=blue]{${0}$}}
      child[color=green] {node [circle,draw,color=green]{${2}$}
        child[color=green] {node[circle,draw,color=green] {${2}$}
          child[color=black] {node[circle,draw,color=red] {${1}$}
            child {node[rectangle, draw, color=blue] {${0}$}}
            child {node[rectangle, draw, color=blue] {${0}$}}
          }
          child[color=black] {node[circle,draw,color=red] {${1}$}
            child {node[rectangle, draw, color=blue] {${0}$}}
            child {node[rectangle, draw, color=blue] {${0}$}}
          }
        }
        child[color=black] {node [rectangle,draw,color=blue]{${0}$}
        }
      }
    }
    child {node[circle,draw,color=green] {${2}$}
      child[color=green] {node [circle,draw,color=green]{${2}$}
        child [color=black]{node[circle,draw,color=red] {${1}$}
          child[color=black] {node[rectangle,draw,color=blue] {${0}$}}
          child[color=red] {node [circle,draw,color=red]{${1}$}
            child[color=black] {node[rectangle,draw,color=blue] {${0}$}}
            child [color=black]{node[rectangle,draw,color=blue] {${0}$}}
          }
        }
        child[color=black] {node [circle,draw,color=red]{${1}$}
          child [color=red]{node [circle,draw,color=red]{${1}$}
            child [color=black]{node [rectangle,draw,color=blue]{${0}$}}
            child [color=black]{node [rectangle,draw,color=blue]{${0}$}}
          }
          child {node[rectangle,draw,color=blue] {${0}$}}
        }
      }
      child {node [circle,draw,color=red]{${1}$}
        child {node[rectangle,draw,color=blue] {${0}$}}
        child {node[rectangle,draw,color=blue] {${0}$}}
      }
    };
  \end{tikzpicture}
  \caption{Binary tree with colored $r$-branches}
  \label{fig:reg-branches}
\end{figure}
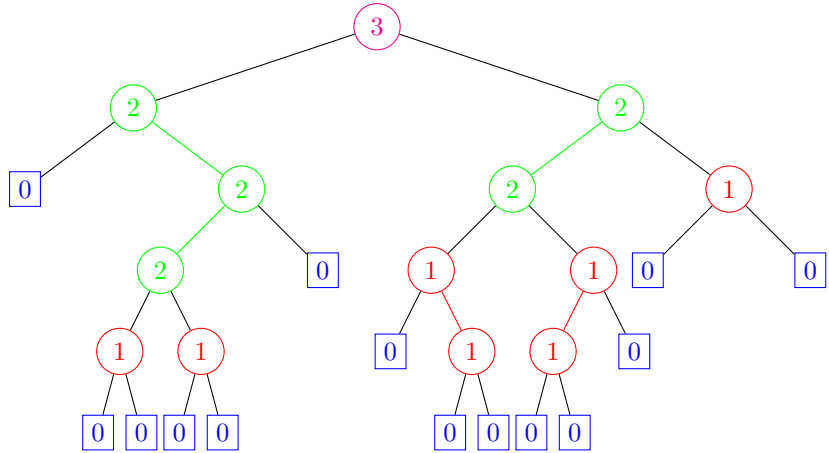

The register function is recursively defined: for the binary tree consisting of
only a leaf we have $\Reg(\square) = 0$, and if a binary
tree $t$ has subtrees $t_{1}$ and $t_{2}$, then the register function is defined to be
\[ \Reg(t) = \begin{cases}
    \max\{\Reg(t_{1}), \Reg(t_{2})\} & \text{ for } \Reg(t_{1}) \neq \Reg(t_{2}),\\
    \Reg(t_{1}) + 1 & \text{ otherwise.}
  \end{cases}\]
In particular, the numbers attached to the nodes in Figures~\ref{fig:trimming} and
\ref{fig:reg-branches} represent the register function of the subtree rooted at the respective
node. 

Historically, the idea of the register function originated (as the Horton-Strahler
numbers) in \cite{Horton45, Strahler52} in the study of the complexity of river
networks. However, the very same concept also occurs within a computer science
context: arithmetic expressions with binary operators can be expressed as a binary tree
with data in the leaves and operators in the internal nodes. Then, the register function of this
binary expression tree corresponds to the minimal number of registers needed to evaluate
the expression.

There are several publications in which the register function and related concepts are
investigated in great detail, for example Flajolet, Raoult, and Vuillemin~\cite{FlRaVu79},
Kemp~\cite{Kemp79}, Flajolet and Prodinger~\cite{FlPr86}, Louchard and
Prodinger~\cite{Louchard-Prodinger:2008:register-lattice}, Drmota and
Prodinger~\cite{DrPr06}, and Viennot~\cite{viennot:2002:strahler-bijection}. For a
detailed survey on the register function and related topics see~\cite{prodinger:register-survey}.

We continue by observing that the compactification $\Phi$ is a very natural operation
regarding the register function:

\begin{proposition}\label{prop:reg-compact}
  Let $t$ be a binary tree with $\Reg(t) = r \geq 1$. Then $\Phi(t)$ is well-defined and
  the register function of the compactified tree is $\Reg(\Phi(t)) = r - 1$.
\end{proposition}

As an immediate consequence of Proposition~\ref{prop:reg-compact} we find that $\Phi$ can
be applied $r$ times repeatedly to some binary tree $t$ if and only if $\Reg(t) \geq r$
holds. In particular, we obtain
\begin{equation}
  \label{eq:reg-characterization}
  \Phi^{r}(t) = \square \quad \iff \quad \Reg(t) = r.
\end{equation}
With~\eqref{eq:reg-characterization}, the link between the generating functions $B_{r}(z)$
from above and the register function becomes clear: $B_{r}(z)$ is exactly the generating
function of binary trees with register function $\leq r$. 

In order to analyze these recursively defined generating functions an explicit
representation is convenient. As it turns out, the substitution $z = \frac{u}{(1 +
  u)^{2}} =: Z(u)$ is a helpful tool in this context.

In particular, it can be shown that applying $z \mapsto \frac{z^{2}}{(1 - 2z)^{2}}$
corresponds to $u\mapsto u^{2}$, which helps to find the explicit representation
\[ B_{r}(z) = \frac{1 - u^{2}}{u} \sum_{j=0}^{r} \frac{u^{2^{j}}}{1 - u^{2^{j+1}}}. \]

Note that at this point, we could obtain the generating function for binary trees with
register function equal to $r$ simply by computing the difference $B_{r}(z) - B_{r-1}(z)$
for $r \geq 1$. These functions can be used to study the asymptotic behavior of the
average register function value.

However, as these results are well-known (cf.~\cite{FlRaVu79}), we will continue in a
different direction by studying the number of so-called $r$-branches.

\subsection{$r$-branches}
\label{sec:r-branches}

The register function associates a value to each node (internal nodes as well as leaves), and the
value at the root is the value of the register function of the tree. An $r$-branch is a
maximal chain of nodes labeled $r$. This must be a chain, since the merging of two such
chains would already result in the higher value $r+1$. The nodes of the tree are
partitioned into such chains, from $r=0,1,\ldots$. The goal of this section is the study
of the parameter ``number of $r$-branches'', in particular, the average number of them,
assuming that all binary trees of size $n$ are equally likely.  

This parameter was the main object of the paper \cite{YaYa10}, and some partial results
were given that we are now going to extend. In contrast to this paper, our approach relies
heavily on generating functions which, besides allowing us to verify the results in a
relatively straightforward way, also enables us to extract explicit formul\ae{} for the
expectation (and, in principle, also for higher moments).

A parameter that was not investigated in \cite{YaYa10} is the total number of
$r$-branches, for any $r$, i.e., the sum over $r\ge0$. Here, asymptotics are trickier, and
the basic approach from \cite{YaYa10} cannot be applied. However, in this paper we use the
Mellin transform, combined with singularity analysis of generating functions, a multi-layer
approach that also allowed one of us several years ago to solve a problem by Yekutieli and
Mandelbrot, cf.~\cite{Prodinger97}. The origins of singularity analysis can be found in
\cite{FlOd90}, and for a detailed survey see \cite{FlSe09}.

For reasons of comparisons, let us mention that the value of register function in
\cite{YaYa10} are one higher than here, and that $n$ generally refers there to the number
of leaves, not nodes as here. 

According to our previous considerations, after $r$ iterations of $\Phi$, the $r$-branches
become leaves (or, equivalently, $0$-branches). The bivariate generating function allowing us to count the
leaves of the binary trees is $vB(zv)$.

The proofs of the statements in this section, together with SageMath worksheets containing
the corresponding computations, can be found in the full version of this paper.

\begin{theorem}\label{thm:r-branches} 
  Let $r\in \N_{0}$ be fixed. The expected number of $r$-branches in binary trees of size
  $n$ and the corresponding variance have the following asymptotic expansions:
  \begin{align}
    \label{eq:asy-exp-r-branch}
    E_{n; r} & = \frac{n}{4^{r}} + \frac{1}{6} \Big(1 +
    \frac{5}{4^{r}}\Big) + \frac{1}{20 n}\Big(4^{r} - \frac{1}{4^{r}}\Big) + \frac{1}{12
    n^{2}} \Big(\frac{5\cdot 16^{r}}{21}  - \frac{7\cdot 4^{r}}{10}  + \frac{97}{210\cdot 4^{r}}\Big)
    + O(n^{-3}), \\
    \label{eq:asy-var-r-branch}
    V_{n;r} & = \frac{4^r-1}{3\cdot 16^r} n - \frac{2\cdot 16^{r} - 25\cdot 4^{r} + 23}{90\cdot
    16^{r}} - \frac{13\cdot 64^{r} - 14\cdot 16^{r} + 7\cdot 4^{r} - 6}{420\cdot 16^{r} n} + O(n^{-2}).
  \end{align}
\end{theorem}
Of course, the expected number of $r$-branches can also be computed explicitly by using
Cauchy's integral formula. This yields the following result:

\begin{proposition}
  \label{prop:explicit-exp-r-branch}
  The expected number of $r$-branches in binary trees of size $n$ is given by the explicit
  formula
  \begin{equation}
    \label{eq:explicit-exp-r-branch}
    E_{n;r} = \frac{n + 1}{\binom{2n}{n}} \sum_{\lambda \geq 1} \lambda
    \bigg[\binom{2n}{n + 1 - \lambda 2^{r}} - 2\binom{2n}{n - \lambda 2^{r}} +
    \binom{2n}{n - 1 - \lambda 2^{r}}\bigg].
  \end{equation}
\end{proposition}

\subsection{The total number of branches}
\label{sec:all-branches}

So far, we were dealing with fixed $r$, and the number of $r$-branches in trees of size
$n$, for large $n$. Now we consider the total number of such branches, i.e., the sum over
$r\ge0$, which was not considered in \cite{YaYa10}. First, to get an explicit formula, the
results from Proposition~\ref{prop:explicit-exp-r-branch} can be summed.

\begin{corollary}\label{cor:explicit-exp-branches} 
  The expected number of branches in binary trees of size $n$, denoted as $E_{n}$, is
  given by the explicit formula
  \begin{equation*}
    E_{n} = \frac{n+1}{\binom{2n}{n}}\sum_{k=1}^{n+1} (2 - 2^{-v_{2}(k)}) k \bigg[\binom{2n}{n+1-k}-2\binom{2n}{n-k}+
    \binom{2n}{n-1-k}\bigg],
  \end{equation*}
  where $v_{2}(k)$ is the dyadic valuation of $k$, i.e., the highest exponent
  $\nu$ such that $2^\nu$ divides $k$.
\end{corollary}
While it is absolutely possible to work out the asymptotic growth from this explicit
formula, at it was done in earlier papers \cite{FlRaVu79, Kemp79}, we choose a faster
method, like in \cite{FlPr86}. It works on the level of generating functions and uses the
Mellin transform together with singularity analysis of generating functions~\cite{FlSe09,
  bona:prodinger:2015:analyt}.

The following theorem describes the asymptotic behavior for the expected number of
branches in a binary tree.

\begin{theorem}
  \label{thm:asy-branches}
  The expected value of the total number of branches in a random binary tree of size $n$
  admits the asymptotic expansion
  \begin{equation*}
    E_{n} = \frac{4n}{3} + \frac{1}{6}\log_{4} n - \frac{2 \zeta'(-1)}{\log2} -
    \frac{\gamma}{12\log 2} - \frac{1}{6\log 2} + \frac{43}{36}
    +\delta(\log_4n)+ O \Big(\frac{\log n}{n}\Big),
  \end{equation*}
  where 
  \begin{equation*}
    \delta(x) := \frac1{\log2}\sum_{k\neq 0}
    \Gamma\Big(\frac{\chi_k}{2}\Big)\zeta(\chi_k-1)(\chi_k-1)e^{2\pi i k x}
  \end{equation*}
  is a $1$-periodic function of mean zero, given by its Fourier series expansion.
\end{theorem}
\begin{remark}
  Note that the value of the derivative of the zeta function is given by $\zeta'(-1) = -
  \frac{1}{12} - \log A \approx -0.1654211437$, where $A$ is the Glaisher-Kinkelin
  constant (cf.~\cite[Section 2.15]{Finch:constants:2003}).
\end{remark}
\newcommand{\figbranchesfluc}{
\begin{figure}[htbp]
  \centering
  \begin{tikzpicture}
    \begin{axis}[xlabel={$x$}, ylabel={$\delta(x)$}, width=14cm, height=8cm,
      yticklabel style={/pgf/number format/fixed, /pgf/number format/precision=3},
      xtick = {2,2.5,...,5}, legend entries = {empirical, Fourier series},
      ytick = {-0.09, -0.06, ..., 0.06}, legend pos = south east]
      \addplot+[only marks, mark size=0.4pt] table[x=x, y=empirical] {theorem2_fluc.dat};
      \addplot+[no marks] table[x=x, y=fourier] {theorem2_fluc.dat};
    \end{axis}
  \end{tikzpicture}
  \caption{Partial Fourier series (20 summands) compared with the empirical values of the
    function $\delta$ from Theorem~\ref{thm:asy-branches}}
  \label{fig:branches-fluc}
\end{figure}
}

\begin{remark}
  The occurrence of the periodic fluctuation $\delta$ where the argument is logarithmic in
  $n$ is actually not surprising: while this phenomenon is already very common in the
  context of the register function, fluctuations appear very often in the asymptotic
  analysis of sums. 
\end{remark}

While this multi-layer approach enabled us to analyze the expected value of the number of
branches in binary trees of size $n$, the same strategy fails for computing the
variance. This is because the random variables modeling the number of $r$-branches are
correlated for different values of $r$---and thus, the sum of the variances (which we
compute by our approach) differs from the variance of the sum.

This concludes our study of the number of branches per binary tree. In the next section,
we analyze a quantity that has similar properties as the register function, but is defined on
simple two-dimensional lattice paths.

\section{A Similar Recursive Scheme Involving Lattice Paths}
\label{sec:paths}

Recall that the register function describes the number of compactifications of a binary
tree required in order to reduce the tree to a leaf. By defining a similar
process for simple two-dimensional lattice paths, a
function that plays a similar role as the register function is obtained.

Simple two-dimensional lattice paths are sequences of the symbols $\{\upa,
\righta, \downa, \lefta\}$. It is easy to see that the generating function
counting these paths (without the path of length $0$) is
\[ L(z) = \frac{4z}{1 - 4z} = 4z + 16z^{2} + 64z^{3} + 256z^{4} + 1024z^{5} + \cdots.  \]

\begin{proposition}\label{prop:lat-relation}
  The generating function $L(z) = \frac{4z}{1-4z}$ fulfills the functional equation
  \begin{equation}\label{eq:lat-relation}
    L(z) = 4 L\Big(\frac{z^{2}}{(1 - 2z)^{2}}\Big) + 4z.
  \end{equation}
\end{proposition}
\begin{remark}
  It is easy to verify this result by means of substitution and expansion. However, we want
  to give a combinatorial proof---this approach also motivates the definition of a
  recursive generation process for lattice paths, similar to the process for binary trees
  from above.
\end{remark}
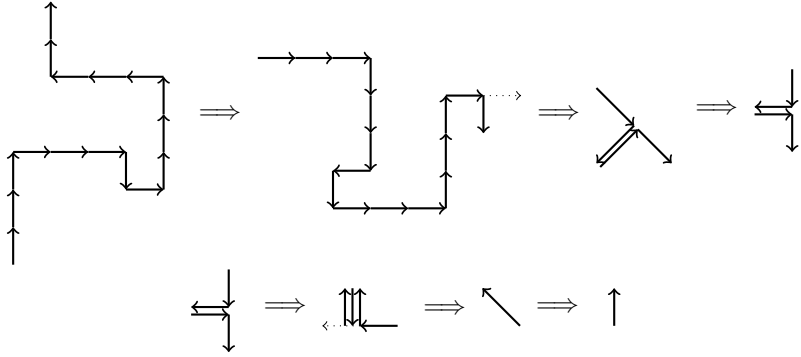
\begin{figure}[htbp]
    \centering
    \begin{tikzpicture}[scale=0.5]
      \foreach \x  [remember=\x as \lastx (initially {(0,0)})] in %
      {(0,1), (0,2), (0,3), (1,3), (2,3), (3,3), (3,2), (4,2), (4,3), (4,4), (4,5), (3,5),
        (2,5), (1,5), (1,6), (1,7)} 
      {
        \draw[->, thick] \lastx -- \x;
      }
      \node at (5.5, 4) {$\Longrightarrow$};
    \end{tikzpicture}\raisebox{0.75cm}{
      \begin{tikzpicture}[scale=0.5]
        \foreach \x  [remember=\x as \lastx (initially {(0,0)})] in %
        {(1,0), (2,0), (3,0), (3, -1), (3, -2), (3, -3), (2, -3), (2, -4), (3,-4), (4, -4),
          (5,-4), (5, -3), (5, -2), (5, -1), (6, -1)} 
        {
          \draw[->, thick] \lastx -- \x;
        }
        \draw[->, dotted] (6, -1) -- (7, -1);
        \draw[->, thick] (6, -1) -- (6, -2);
        \node at (8, -1.5) {$\Longrightarrow$};
      \end{tikzpicture}}\raisebox{1.3cm}{
      \begin{tikzpicture}[scale=0.5]
        \draw[->, thick] (0,0) -- (1,-1);
        \draw[->, thick] (1, -1) to (0, -2);
        \draw[->, thick] (0.1, -2.1) to (1.1, -1.1);
        \draw[->, thick] (1.1, -1.1) -- (2, -2);
        \node at (3.2, -0.55) {$\Longrightarrow$};
      \end{tikzpicture}}\raisebox{1.5cm}{
      \begin{tikzpicture}[scale=0.5]
        \draw[->, thick] (0,0) -- (0,-1);
        \draw[->, thick] (0, -1) to (-1, -1);
        \draw[->, thick] (-1, -1.2) to (0, -1.2);
        \draw[->, thick] (0, -1.2) -- (0, -2.2);
      \end{tikzpicture}}
    
    \begin{tikzpicture}[scale=0.5]
      \draw[->, thick] (0,0) -- (0,-1);
        \draw[->, thick] (0, -1) to (-1, -1);
        \draw[->, thick] (-1, -1.2) to (0, -1.2);
        \draw[->, thick] (0, -1.2) -- (0, -2.2);
      \node at (1.5, -1) {$\Longrightarrow$};
    \end{tikzpicture}\raisebox{0.35cm}{
      \begin{tikzpicture}[scale=0.5]
        \draw[->, thick] (0,0) -- (-1,0);
        \draw[->, thick] (-1, 0) to (-1, 1);
        \draw[->, thick] (-1.2, 1) to (-1.2, 0);
        \draw[->, dotted] (-1, 0) -- (-2, 0);
        \draw[->, thick] (-1.4, 0) to (-1.4, 1);
        \node at (1.25, 0.45) {$\Longrightarrow$};
      \end{tikzpicture}}\raisebox{0.35cm}{
      \begin{tikzpicture}[scale=0.5]
        \draw[->, thick] (0,0) -- (-1,1);
        \node at (1, 0.5) {$\Longrightarrow$};
      \end{tikzpicture}}\raisebox{0.35cm}{\hspace{0.25cm}
      \begin{tikzpicture}[scale=0.5]
        \draw[->, thick] (0,0) -- (0,1);
      \end{tikzpicture}}
    \caption{Repeated application of the reduction $\Phi_{L}$ on a path with compactification degree 2}
    \label{fig:lat-red-vis}
  \end{figure}
\begin{proof}[(Sketch)]
  While we leave the detailed proof to the full version of this paper, we still want to
  introduce a lattice path reduction which plays an analogous role as the binary tree
  reduction in the proof of Proposition~\ref{prop:tree-identity}.

  We consider the reduction $\Phi_{L}$, which acts on any given lattice path $\ell$ with
  length $\geq 2$ as follows:

  First, the path needs to be modified such that it starts horizontally and ends
  vertically. This is done by rotation to the right of the entire path and/or the very
  last step, respectively.

  Then, the horizontally starting and vertically ending path is reduced by replacing
  each pair of successive horizontal-vertical path segments in the following way:
  \begin{itemize}
  \item If a segment starts with $\righta$ and the first vertical step is $\upa$, replace it by
    $\nearrow$,
  \item if a segment starts with $\righta$ and the first vertical step is $\downa$, replace it
    by $\searrow$,
  \item if a segment starts with $\lefta$ and the first vertical step is $\downa$, replace it by
    $\swarrow$,
  \item and if a segment starts with $\lefta$ and the first vertical step is $\upa$, replace it
    by $\nwarrow$.
  \end{itemize}
  Rotating the resulting path by $45^{\circ}$ in order to obtain a path with horizontal
  and vertical steps then yields $\Phi_{L}(\ell)$. It can be shown that this reduction
  corresponds to the right-hand side of~\eqref{eq:lat-relation}.
\end{proof}

The process described in the proof of Proposition~\ref{prop:lat-relation} allows us to
assign a unique number to each lattice path:
\begin{definition}
  Let $\ell$ be a simple two-dimensional lattice path consisting of at least one step. We
  define the \emph{compactification degree} of $\ell$, denoted as $\cdeg(\ell)$ as
  \[ \cdeg(\ell) = n \quad \iff\quad \Phi_{L}^{n}(\ell) \in \{\upa, \righta,
    \downa, \lefta\}.  \]
\end{definition}
\begin{remark}
  The parallels between the compactification degree and the register function are obvious:
  both count the number of times some given mathematical object can be reduced according
  to some rules until an atomic form of the respective object is obtained. Therefore, both
  functions describe, in some sense, the complexity of a given structure.
\end{remark}

In the remainder of this section we want to derive some asymptotic results for the
compactification degree, namely the expected degree of a lattice path of given length as
well as the corresponding variance.

Analogously to our strategy for~\eqref{eq:u1}, we want to interpret~\eqref{eq:lat-relation} as
a recursive generation process as well and therefore set
\begin{equation*}
  \label{eq:lat-recursion}
  L_{0}^{=}(z) = 4z,\quad L_{r}^{=}(z) = 4 L_{r-1}^{=}\Big(\frac{z^{2}}{(1 - 2z)^{2}}\Big),\quad
  r\geq 1.
\end{equation*}
With the help of the substitution $z = Z(u)$ the generating function can be written
explicitly as
\begin{equation}
  \label{eq:lat-gen-r}
  L_{r}^{=}(z) = 4^{r+1} \frac{u^{2^{r}}}{(1 + u^{2^{r}})^{2}}.
\end{equation}

The coefficients of this function can be extracted explicitly by applying Cauchy's
integral formula.
\begin{proposition}\label{prop:lat-coef-r}
  The number of two-dimensional simple lattice paths of length $n$ that have
  compactification degree $r$ is given by
  \[ [z^{n}] L_{r}^{=}(z) = 4^{r+1} \sum_{\lambda \geq 0} \lambda (-1)^{\lambda - 1}
    \bigg[\binom{2n-1}{n - \lambda 2^{r}} - \binom{2n-1}{n - \lambda 2^{r} - 1}\bigg].  \]
\end{proposition}

In fact, by studying the substitution $z = Z(u)$ closely, the asymptotic behavior of
the coefficients of $L_{r}^{=}(z)$ can be extracted as well.

We turn to the investigation of the expected
compactification degree. Let $\mathcal{L}_{n}$ denote the set of simple two-dimensional lattice
paths of size $n$. Consider the family of random variables $X_{n}\colon \mathcal{L}_{n} \to
\mathbb{N}_{0}$ modeling the compactification degree of the lattice paths of length $n$
under the assumption that all paths are equally likely. The following results are
immediate consequences of Proposition~\ref{prop:lat-coef-r}.

\begin{corollary}
  \label{cor:lat-random-explicit}
  The probability that a lattice path of length $n$ has compactification degree $r$ is
  given by the explicit formula
  \[ \P(X_{n} = r) = \frac{[z^{n}] L_{r}^{=}(z)}{4^{n}} = 4^{r+1-n} \sum_{\lambda \geq 0}
    \lambda (-1)^{\lambda - 1} \bigg[\binom{2n-1}{n - \lambda 2^{r}} - \binom{2n-1}{n -
      \lambda2^{r} -1}\bigg],  \]
  and the expected compactification degree for paths of length $n$ is given by
  \begin{equation}\label{eq:lat-exp-explicit}
    \E X_{n} = \sum_{k \geq 1} 8k(2^{v_{2}(k)} - 1) \bigg[\binom{2n-1}{n -
    k} - \binom{2n-1}{n - k - 1}\bigg].
  \end{equation}
\end{corollary}
\begin{remark}
  The formula for $\P(X_{n} = r)$ is very similar to the results for the classical
  register function obtained by Flajolet (cf.~\cite{these-d'etat}). It is likely that
  applying the techniques that were used
  in~\cite{Louchard-Prodinger:2008:register-lattice} could be used to determine expansions
  for arbitrary moments.
\end{remark}
The following theorem characterizes the asymptotic behavior of the expected
compactification degree and the corresponding variance.

\begin{theorem}
  \label{thm:lat-results}
  The expected compactification degree of simple two-dimensional lattice paths of length
  $n$ admits the asymptotic expansion
  \begin{equation}\label{eq:lat-exp-expansion}
  \E X_{n} = \log_{4} n + \frac{\gamma + 2 - 3\log 2}{2 \log 2} +
    \delta_{1}(\log_{4} n) + O(n^{-1}),
  \end{equation}
  and for the corresponding variance we have
  \begin{multline}\label{eq:lat-var-expansion}
    \V X_{n} = \frac{\pi^{2} - 24\log^{2}\pi - 48 \zeta''(0) - 24}{24 \log^{2} 2} -
    \frac{2 \log\pi}{\log 2} - \frac{11}{12} + \delta_{2}(\log_{4}n) \\ - \frac{\gamma + 2 -
      3\log 2}{\log 2} \delta_{1}(\log_{4} n) + \delta_{1}^{2}(\log_{4} n) +
    O\Big(\frac{1}{\log n}\Big) 
  \end{multline}
  where $\delta_{1}(x)$ and $\delta_{2}(x)$ are $1$-periodic fluctuations of mean zero
  whose Fourier coefficients can be given explicitly.
\end{theorem}

\subsection{Fringes}\label{sec:fringes}

We define the $r$th \emph{fringe} of a given lattice path $\ell$ of length
$\geq 1$ to be $\Phi_{L}^{r}(\ell)$, i.e.\ the $r$th fringe is given by the $r$th
reduction of the path. In particular, if $\ell$ can be reduced $r$ times, we call the
length of $\Phi_{L}^{r}(\ell)$ the size of the $r$th fringe. Otherwise, we say that this
size is $0$.

The $r$th fringes of positive size can then be enumerated by the bivariate generating function
\begin{equation*}
  \FGF_r(z, v)=\sum_{\substack{\ell\text{ path} \\\cdeg(\ell)\ge r}} v^{\abs{\Phi_L^r(\ell)}}z^{\abs{\ell}}
\end{equation*}
where $\abs{\ell}$ denotes the length of a lattice path.

It can be shown that $H_{r}(z,v)$ fulfills the recursion
\[ H_{0}(z, v) = \frac{4zv}{1 - 4zv},\quad H_{r}(z, v) = 4H_{r-1}\Big(\Big(\frac{z}{1 -
    2z}\Big)^{2}, v\Big),\ r\geq 1,  \]
which can be used to find the explicit representation
\[ H_{r}(z, v) = \frac{4^{r+1} u^{2^{r}} v}{(1 + u^{2^{r}})^{2} - 4u^{2^{r}} v}.  \]

The generating function $H_{r}(z,v)$ can now be used to derive the asymptotic behavior of
the expectation $E_{n;r}^{L}$ and the variance $V_{n;r}^{L}$ of the size of the $r$th
fringe, where all paths of length $n$ arise with the same probability.

\begin{theorem}\label{thm:fringe-sizes}
  Let $r\in \N_{0}$ be fixed. The expectation and variance of the $r$th fringe size of a
  random path of length $n$ have the asymptotic expansions
  \begin{equation}
    \label{eq:fringes-exp}
    E_{n;r}^{L} = \frac{n}{4^{r}} + \frac{1 - 4^{-r}}{3} + O(n^{3} \theta_{r}^{-n})
  \end{equation}
  and
  \begin{equation}
    \label{eq:fringes-var}
    V_{n;r}^{L} = \frac{4^{r} - 1}{3\cdot 16^{r}} n + \frac{-2\cdot 16^{r} - 5\cdot 4^{r}
      + 7}{45\cdot 16^{r}} + O(n^{5} \theta_{r}^{-n}),
  \end{equation}
  where $\theta_{r} = \frac{4}{2 + 2\cos(2\pi/2^{r})} > 1$.
  If additionally $r > 0$, then for the random variables $Y_{n;r}$
  modeling the $r$th fringe size of lattice paths of length $n$ we have
  \[ \P\bigg(\frac{Y_{n;r} - E_{n;r}}{\sqrt{V_{n;r}}} \leq x\bigg) = \frac{1}{\sqrt{2\pi}}
    \int_{-\infty}^{x} e^{-w^{2}/2}~dw + O(n^{-1/2}),  \]
  i.e.\ the random variables $Y_{n;r}$ are asymptotically normally distributed.
\end{theorem}

As we have the generating function $H_{r}(z,v)$ in an explicit form, the expected value
can also be extracted explicitly by means of Cauchy's integral formula.

\begin{proposition}\label{prop:fringes-explicit}
  For given $r\in \N_{0}$, the $r$th expected fringe size of a random path of
  length $n$ is given by the explicit formula
  \[ E_{n;r}^{L} = 4^{r+1-n} \sum_{\lambda \geq 1} \frac{2\lambda^{3} + \lambda}{3}
    \bigg[\binom{2n-1}{n-2^{r}\lambda} - \binom{2n-1}{n - 2^{r}\lambda -1}\bigg].  \]
\end{proposition}

Analogously to our investigations concerning branches in binary trees, we also
study the asymptotic behavior of the expected fringe size, i.e.\ the sum over the size of
the $r$th fringes for $r \geq 0$. Like the compactification degree, this parameter can
also be interpreted as a complexity measure for lattice paths.

\begin{corollary}\label{cor:fringe-all-explicit}
  The expected fringe size $E_{n}^{L}$ of a random path of length $n$
  can be computed as
  \[ E_{n}^{L} = \frac{1}{12\cdot 4^{n}} \sum_{k=1}^{n} \big(2k^{3} (2 - 2^{-v_{2}(k)}) +
    k(2^{v_{2}(k) + 1} - 1)\big)\bigg[\binom{2n-1}{n-k} - \binom{2n-1}{n-k-1}\bigg].  \]
\end{corollary}

The following theorem quantifies the asymptotic behavior of $E_{n}^{L} := \sum_{r\geq 0} E_{n;r}^{L}$.

\begin{theorem}\label{thm:fringe-size}
  Asymptotically, the behavior of the expected fringe size $E_{n}^{L}$ for a random path
  of length $n$ is given by
  \begin{equation}
    \label{eq:fringe-size-exp}
    E_{n}^{L} = \frac{4}{3} n + \frac{1}{3} \log_{4}n + \frac{5 + 3\gamma - 11\log 2}{18
      \log 2} + \delta(\log_{4} n) + O\Big(\frac{\log n}{n}\Big),
  \end{equation}
  where $\delta(x)$ is a $1$-periodic fluctuation of mean zero with Fourier series
  expansion
  \[ \delta(x) = \sum_{k\neq 0}\frac{2}{3\sqrt{\pi} \log 2} \Gamma\Big(\frac{3 + \chi_{k}}{2}\Big)
    \bigl(2\zeta(\chi_{k} - 1) + \zeta(\chi_{k} + 1)\bigr) e^{2k\pi i x}.  \]
\end{theorem}

\newcommand{\MR}{}

\bibliographystyle{amsplaininitialsurl}
\bibliography{register-reduction,bib/cheub}

\providecommand{\Submitted}{Submitted} \providecommand{\availableat}{ available
  at } \providecommand{\alsoavailableat}{ also available at }
  \providecommand{\evavailableat}{earlier version available at }
  \providecommand{\toappearin}{To appear in } \providecommand{\toappear}{to
  appear} \providecommand{\inpreparation}{in preparation}
  \providecommand{\doi}[1]{\href{http://dx.doi.org/#1}{\path{doi:#1}}}
  \providecommand{\etc}{\emph{etc.}}\def\cprime{$'$}
\providecommand{\bysame}{\leavevmode\hbox to3em{\hrulefill}\thinspace}
\providecommand{\MR}{\relax\ifhmode\unskip\space\fi MR }
\providecommand{\MRhref}[2]{%
  \href{http://www.ams.org/mathscinet-getitem?mr=#1}{#2}
}
\providecommand{\href}[2]{#2}
\begin{thebibliography}{10}

\bibitem{DrPr06}
M.~Drmota and H.~Prodinger, \emph{The register function for {$t$}-ary trees},
  ACM Trans. Algorithms \textbf{2} (2006), no.~3, 318--334.

\bibitem{Finch:constants:2003}
S.~R. Finch, \href{http://dx.doi.org/10.1017/CBO9780511550447
  }{\emph{Mathematical constants}}, Encyclopedia of Mathematics and its
  Applications, vol.~94, Cambridge University Press, Cambridge, 2003.
  \MR{2003519 (2004i:00001)}

\bibitem{these-d'etat}
P.~Flajolet, \emph{Analyse d'algorithmes de manipulation d'arbres et de
  fichiers}, Cahiers du Bureau Universitaire de Recherche Opérationnelle
  \textbf{34/35} (1981), 1--209.

\bibitem{FlPr86}
P.~Flajolet and H.~Prodinger, \emph{Register allocation for unary-binary
  trees}, SIAM J. Comput. \textbf{15} (1986), 629--640.

\bibitem{FlRaVu79}
P.~Flajolet, J.-C. Raoult, and J.~Vuillemin, \emph{The number of registers
  required for evaluating arithmetic expressions}, Theoret. Comput. Sci.
  \textbf{9} (1979), no.~1, 99--125.

\bibitem{FlOd90}
P.~Flajolet and A.~Odlyzko, \emph{Singularity analysis of generating
  functions}, SIAM J. Discrete Math. \textbf{3} (1990), 216--240.

\bibitem{FlSe09}
P.~Flajolet and R.~Sedgewick, \emph{Analytic combinatorics}, Cambridge
  University Press, Cambridge, 2009.

\bibitem{Horton45}
R.~E. Horton, \emph{Erosioned development of systems and their drainage
  basins}, Geol. Soc. Am. Bull. \textbf{56} (1945), 275--370.

\bibitem{Kemp79}
R.~Kemp, \emph{The average number of registers needed to evaluate a binary tree
  optimally}, Acta Inform. \textbf{11} (1978/79), no.~4, 363--372.

\bibitem{Louchard-Prodinger:2008:register-lattice}
G.~Louchard and H.~Prodinger, \emph{The register function for lattice paths},
  Fifth {C}olloquium on {M}athematics and {C}omputer {S}cience, Discrete Math.
  Theor. Comput. Sci. Proc., AI, Assoc. Discrete Math. Theor. Comput. Sci.,
  Nancy, 2008, pp.~135--148. \MR{2508783}

\bibitem{Prodinger97}
H.~Prodinger, \emph{On a problem of {Y}ekutieli and {M}andelbrot about the
  bifurcation ratio of binary trees}, Theoret. Comput. Sci. \textbf{181}
  (1997), no.~1, 181--194.

\bibitem{bona:prodinger:2015:analyt}
H.~Prodinger, \href{https://www.crcpress.com/product/isbn/9781482220858
  }{\emph{Analytic methods}}, Handbook of Enumerative Combinatorics
  (M.~B{\'o}na, ed.), CRC Press Series on Discrete Mathematics and its
  Applications, Chapman \& Hall/CRC, Boca Raton, FL, 2015, pp.~173--252.

\bibitem{prodinger:register-survey}
\bysame,
  \href{http://math.sun.ac.za/~hproding/pdffiles/flajolet-introduction-register.pdf
  }{\emph{Introduction to {P}hilippe {F}lajolet's work on the register function
  and related topics}}, Philippe Flajolet's Collected Papers (M.~D. Ward, ed.),
  vol.~V, to appear.

\bibitem{shapiro-1976}
L.~W. Shapiro, \emph{A short proof of an identity of {T}ouchard's concerning
  {C}atalan numbers}, J. Combinatorial Theory Ser. A \textbf{20} (1976), no.~3,
  375--376. \MR{0406819 (53 \#10605)}

\bibitem{Strahler52}
A.~N. Strahler, \emph{Hypsomic analysis of erosional topography}, Geol. Soc.
  Am. Bull. \textbf{63} (1952), 1117--1142.

\bibitem{sage-7.0}
{The Sage Developers}, \emph{{S}age {M}athematics {S}oftware ({V}ersion 7.0)},
  2016, {\tt http://www.sagemath.org}.

\bibitem{touchard-1924}
J.~Touchard, \emph{Sur certaines \'equations fonctionnelles}, Proc. Internat.
  Math. Congress, vol. I (1928), 1924, pp.~465--472.

\bibitem{viennot:2002:strahler-bijection}
X.~G. Viennot, \href{http://dx.doi.org/10.1016/S0012-365X(01)00265-5 }{\emph{A
  {S}trahler bijection between {D}yck paths and planar trees}}, Discrete Math.
  \textbf{246} (2002), no.~1-3, 317--329, Formal power series and algebraic
  combinatorics (Barcelona, 1999). \MR{1887493}

\bibitem{YaYa10}
K.~Yamamoto and Y.~Yamazaki, \emph{Topological self-similarity on the random
  binary-tree model}, J. Stat. Phys. \textbf{139} (2010), no.~1, 62--71.

\end{thebibliography}

\end{document}